\documentclass[12pt]{amsart}
\usepackage{amsmath, amssymb}

\newtheorem{theorem}{Theorem}
\newtheorem{corollary}{Corollary}
\newtheorem{lemma}{Lemma}
\newtheorem{proposition}{Proposition}

\theoremstyle{definition}

\newtheorem{example}{Example}

\numberwithin{theorem}{section}
\numberwithin{corollary}{section}
\numberwithin{lemma}{section}
\numberwithin{proposition}{section}
\numberwithin{definition}{section}
\numberwithin{example}{section}

\newcommand{\p}{\mathbf{P}} % bold P
\newcommand{\aff}{\mathbf{A}} % bold A
\newcommand{\calO}{{\mathcal O}}

\DeclareMathOperator{\Spec}{Spec}

\DeclareMathOperator{\Proj}{Proj}
\DeclareMathOperator{\height}{height}

\begin{document}

\title[Flat Normal Cones and Segre Classes]
    {Detecting Flat Normal Cones \\ Using Segre Classes}

\author[S.\ J.\ Colley]{Susan Jane Colley}
\address{Department of Mathematics,
Oberlin College, Oberlin, Ohio 44074, USA}
\email{sjcolley@math.oberlin.edu}
\thanks{Research of the first author was partly supported
by National Science Foundation grant \#DMS-9805819.}

\author[G.\ Kennedy]{Gary Kennedy}
\address{Ohio State University at Mansfield, 1680 University Drive,
Mansfield, Ohio 44906, USA}
\email{kennedy@math.ohio-state.edu}



\begin{abstract}
    Given a flat, projective morphism $Y \to T$ from an
    equidimensional scheme to a nonsingular curve and a
    subscheme $Z$ of $Y$, we give conditions under which
    specialization of the Segre class $s(N_{Z}Y)$ of the normal cone
    of $Z$ in $Y$ implies flatness of the normal cone.  We apply this
    result to study when the relative tangent star cone of a flat
    family is flat.
\end{abstract}

\maketitle

%
%
%% Main body of paper starts here %%
%
%

\section{Introduction} \label{intro}

This paper, situated on the cusp of commutative algebra and algebraic
geometry, is concerned with the geometry of families of blowups.
Recall that the \emph{blowup} of an affine scheme
$\Spec{A}$ along a subscheme
$\Spec{A/I}$ is $\Proj(R_I(A))$, where $R_I(A)$ is the
\emph{Rees algebra} $A \oplus I \oplus I^2 \oplus \cdots$,
and that its \emph{exceptional divisor}
is $\Proj(\operatorname{gr}_I(A))$, where
$\operatorname{gr}_I(A)$ is the \emph{associated graded ring}
$A/I \oplus I/I^2 \oplus \cdots$.
The exceptional divisor can also be regarded as the
projectivization of the \emph{normal cone}
$\Spec(\operatorname{gr}_I(A))$.
In algebraic geometry the \emph{center}
$\Spec{A/I}$ is often a nonsingular
variety, but the definitions make sense
in general.
They also make sense for an arbitrary scheme $X$,
if one uses the structure sheaf $\calO_{X}$
in place of the ring $A$ and if $I$ is now
regarded as the ideal sheaf of a subscheme.

\par
Using Segre classes, we will examine how the
normal cone varies in a one-parameter family.
Suppose that $Y \to T$ is a flat,
projective morphism from an equidimensional scheme to
a nonsingular curve; let $Z$ be a subscheme of $Y$.
We show that, under a certain hypothesis, if the Segre class
of the normal cone of $Z$ in $Y$ specializes to the corresponding
class for the fiber $Z_t$ in $Y_t$ over a point of $t \in T$,
then the normal
cone is flat above this point. (See Theorem~\ref{mainthm}
for a precise statement.)

\par
Modern intersection theory (as developed by Fulton and
MacPherson in \cite{fulton}) makes Segre classes a central notion.
These classes play the same role for cone
bundles that Chern classes do for vector bundles: they measure twisting
over the base.  Moreover, Segre classes are sensitive to flatness. If
$C$ is a bundle of cones over the $T$-scheme $X$ and if $C$ is flat
over $T$ above a point $t$, then the Segre class of $C$ over $X$ (an
element of $A_* X$) specializes to the Segre class of $C_t$ over
$X_t$. In a very rough sense, our Theorem~\ref{mainthm} reverses
this implication in the case of normal cones.

\par
Because we hope that our results will be of interest to both algebraic
geometers and commutative algebraists, we need to review some
fundamental notions.  Thus in {\S}\ref{segre} we briefly recall basic
definitions about cones and their Segre classes.
In {\S}\ref{intflat} we describe our algebraic notion
of ``internal flatness'' (which follows either from flatness or from
Serre's $S_1$ condition) and provide a computational criterion for
detecting it.  In {\S}\ref{detectflat} we present our main theorem.
Finally, in {\S}\ref{examples},
we examine various examples in the case where $Y$ is
a fiber product $X\times_T X$ and the subscheme $Z$ is $X$
embedded as the diagonal.  The normal cone in this case is
known as the \emph{(relative) tangent star cone}.
These examples suggest that
Theorem~\ref{mainthm} is probably optimal.

\par
To be fully useful, our Theorem~\ref{mainthm} requires a
``front end'' guaranteeing that certain families of normal
cones satisfy the internal flatness condition.
We know of one result in this direction: the tangent star cone
of a hypersurface is Cohen-Macaulay \cite{gkflat}.
Other such results would be desirable.
On the other hand, for tangent star cones the Segre
class specialization condition is readily checked,
since---as explained at the beginning of {\S}4 of
\cite{ky}---the Segre classes of a projective scheme can be computed
from the double point classes of generic linear projections.

\par
We are pleased to thank Mel Hochster, Claudia Miller, Bernd Ulrich, 
and Wolmer Vasconcelos for their helpful communications with us.

\section{Cones and Segre classes} \label{segre}

\par
All schemes considered in this paper are over a field of 
characteristic zero.  If $S$ is a graded sheaf of $\calO_{X}$-algebras
on a scheme $X$, we say that $C=\Spec(S)$ is a \emph{cone}, and we call
the schemes $\p(C)=\Proj(S)$ and $\p(C\oplus1)=\Proj(S\oplus\mathcal{O}_{X})$
its \emph{projectivization} and \emph{projective completion},
respectively.
Following Chapter 4 of \cite{fulton} or {\S}2.3 of
\cite{joins}, we define the
\emph{Segre class} of $C$ to be
\begin{equation*}
s(C) = \sum_{i \geq 0}q_*(c_1\left(\mathcal{O}(1))^i \cap [\p(C\oplus1)]\right),
\end{equation*}
an element of the Chow group of algebraic cycles on $X$. (Here $q$
denotes the projection map from $\p(C\oplus1)$ to $X$.)

\par
Among the basic properties of the Segre class, we wish to take special
note of its behavior under specialization (Example 10.1.10 of
\cite{fulton}): if $X$ is a scheme over $T$, with fiber $X_t$ over a
particular point $t$, and if $C$ is flat over $T$, then $s(C)$
specializes to $s(C_t)$. Our Theorem~\ref{mainthm} can perhaps
be regarded as a kind of converse, in the case of normal cones.

\section{Internal flatness} \label{intflat}

If $(A,\mathfrak{m})$ is a discrete valuation ring,
we call $T=\Spec(A)$  a \emph{nonsingular curve germ};
its \emph{closed point} is $t=\Spec(A/\mathfrak{m}A)$.
For a morphism $Y \to T$, the fiber over $t$
will be denoted by $Y_t$. If $Y$ is affine, say $Y=\Spec(B)$,
we say that the morphism is \emph{internally flat} if each embedded
prime of $B$ contracts to zero in $A$. To say this geometrically,
the morphism is internally flat if each embedded component of $Y$
dominates $T$. More generally (for $Y$ not necessarily affine),
we say that a morphism $Y \to T$ is
internally flat if for each affine open subset of $Y$ the induced
morphism to $T$ is internally flat.
Note that if $Y$ has no embedded components (in particular, if it is
Cohen-Macaulay) or if the morphism $Y \to T$ is flat, then it
is internally flat. Here is a computational criterion:

\par
\begin{proposition} \label{intflatprop}
Suppose $Y \to T$ is a morphism of affine varieties corresponding
to the ring homomorphism $A \to B$, where $A$ is a discrete valuation
ring with uniformizing parameter $\tau$
and $B=A[x_1,\dots,x_n]/I$ is equidimensional. Let $J \subset I$
be an ideal generated by a regular sequence whose length is
$\height(I)$. Then
\begin{enumerate}
   \item 
   $Y$ has no embedded components
   $\iff J \colon (J \colon I) \subset I$;
   
   \item 
   $Y \to T$ is flat $\iff I \colon \tau \subset I$;
   
   \item 
   $Y \to T$ is internally flat 
   $\iff (I \colon \tau) \cap (J \colon (J \colon I)) \subset I$.
\end{enumerate}
\end{proposition}
\begin{proof}
    Given an irredundant primary decomposition of $I$,
    the ideal $J \colon (J \colon I)$ is
    the intersection of the
    minimal primary components. (See Proposition 3.3.1 of
    \cite{methods}; cf.\ {\S}2.1 of \cite{construct}.)
    The ideal $I \colon \tau$
    is the intersection of the primary components contracting to zero.
    Thus $(I \colon \tau) \cap (J \colon (J \colon I))$ is the
    intersection of the minimal primary components
    and those contracting to zero.
\end{proof}

\par
The following proposition 
is similar to the Corollary to Theorem 23.1 of \cite{matsu}.

\par
\begin{proposition} \label{ifprop}
Suppose that  $f\colon Y \to T$
is an internally flat morphism from an equidimensional scheme to a
nonsingular curve germ.  Suppose that the fiber $Y_t$ has the expected
dimension, namely, $\dim{Y} - 1$.  Then $f$ is flat.
\end{proposition}

\begin{proof}
We may assume that $Y$ is affine.
Thus $f\colon Y \to T$
corresponds to a homomorphism
$\varphi\colon A \to B$ from a discrete valuation ring
to an equidimensional ring.
Let $F=B \otimes_A A/\mathfrak{m}A$ be the fiber
ring, and let $\tau \in A$ be a uniformizing parameter.
If $\varphi$ is not flat, then $\varphi(\tau)$ is a zero divisor
and thus is an element of some associated prime $\mathfrak{q}$ of $B$.
Since $\mathfrak{q}$ does not contract to zero, it must be minimal.
Thus $\dim{F} \geq \dim{(B/\mathfrak{q})} = \dim{B}$.
\end{proof}

\section{Detecting flatness} \label{detectflat}

\par
Now we state our main theorem.

\begin{theorem} \label{mainthm}

Suppose that:
\begin{enumerate}
\item 
$Y \to T$ is a flat, projective morphism from an 
equidimensional scheme to a nonsingular curve germ (with
closed point $t$) and $Z$ is a subscheme of $Y$;

\item 
the normal cone $N_{Z}Y$ is internally flat over $T$;

\item 
the Segre class $s(N_{Z}Y)$ specializes to $s(N_{Z_t}Y_t)$.
\end{enumerate}
\noindent
Then $N_{Z}Y$ is flat over $T$, and $(N_{Z}Y)_t = N_{Z_t}Y_t$.
\end{theorem}

\par
We first prove a weak version of Theorem~\ref{mainthm},
with an additional hypothesis and a weaker conclusion.

\begin{lemma} \label{bootstrap}
Suppose that:
\begin{enumerate}
   \item 
   $Y \to T$ is a flat, projective morphism from an
   equidimensional scheme to a nonsingular curve germ
   and $Z$ is a subscheme of $Y$;
	 
   \item 
   $N_{Z}Y$ is internally flat over $T$;
   
   \item 
   $s(N_{Z}Y)$ specializes to $s(N_{Z_t}Y_t)$;
   
   \item 
   $Z$ is nowhere dense in $Y$.
\end{enumerate}
\noindent
Then $\p N_{Z}Y$ is flat over $T$, and
$(\p N_{Z}Y)_t = \p N_{Z_t}Y_t$.
\end{lemma}

\begin{proof}
Assume hypotheses 1, 2, 4 of Lemma~\ref{bootstrap}, but that
$\p N_{Z}Y$ is not flat over $T$.  We will first show that
$\p N_{Z}Y$ does not specialize to $\p N_{Z_t}Y_t$, and subsequently
that the Segre class fails to specialize (i.e., that hypothesis 3 does
not hold).

\par
Since $Y$ is equidimensional, so is $\p N_{Z}Y$.
And since $N_{Z}Y \to T$ is internally flat, so is
the morphism $\p N_{Z}Y \to T$.
Thus by Proposition~\ref{ifprop},
$$
\dim (\p N_{Z}Y)_t > \dim \p N_{Z}Y - 1,
$$
and thus the dimension of $(\p N_{Z}Y)_t$ exceeds that
of $\p N_{Z_t}Y_t$.  This tells us that
$\p N_{Z_t}Y_t$ is properly contained
in $(\p N_{Z}Y)_t$.

\par
Now we turn to Segre classes.  In view of hypothesis 4 and 
Example 4.1.2 of \cite{fulton}, the Segre class of the 
normal cone may be calculated using
\begin{equation}\label{segreeqn}
s(N_{Z}Y) =
\sum_{i \geq 0}p_*\left(c_1(\mathcal{O}(1))^i \cap [\p N_{Z}Y]\right)
\end{equation}
where $p\colon \p N_{Z}Y \to Z$ denotes projection.
Because $\p N_{Z_t}Y_t$ is a proper subscheme of
$(\p N_{Z}Y)_t$, the blowup $Bl_{Z_t}Y_t$
is likewise a proper subscheme of $(Bl_{Z}Y)_t$.
Both of these schemes are equidimensional and have the
same dimension, namely, $\dim Y - 1$.
Hence the difference of fundamental classes
$$
D=\left[(Bl_{Z}Y)_t\right] \, - \, \left[Bl_{Z_t}Y_t \right]
$$
must be a positive cycle of dimension $\dim Y - 1$.  Moreover, since
the two cycles comprising $D$ can only disagree over the center of the
blowup, $D$ must be supported on $(\p N_{Z}Y)_t$.

\par
The exceptional (Cartier) divisor $\mathcal{E}$
on $Bl_{Z}Y$ restricts to $c_1(\mathcal{O}(-1))$ on $\p N_{Z}Y$.
Thus
$$
-[(\p N_{Z}Y)_t] + [\p N_{Z_t}Y_t] = -\mathcal{E} \cap D =
   c_1(\mathcal{O}(1)) \cap D.
$$
(Note that the first cap product takes place on the blowup, and
the second in $(\p N_{Z}Y)_t$.)  Thus, by (\ref{segreeqn}), 
the difference between the Segre classes of 
$N_{Z_t}Y_t$ and $(N_{Z}Y)_t$ is given by
$$
s(N_{Z_t}Y_t) - s((N_{Z}Y)_t) =
\sum_{i \geq 0}\pi_*\left(c_1(\mathcal{O}(1))^{i+1} \cap D\right),
$$
where $\pi\colon (\p N_{Z}Y)_t \to Z_t$ is projection.
To show that this cycle class is nonzero, we consider an
arbitrary component $V$ of $D$.  Set $d = \dim V - \dim \pi(V)$.
Then
$$
\pi_*\left(c_1(\mathcal{O}(1))^j \cap V\right) = 0
$$
if $j < d$, while
$$
\pi_*\left(c_1(\mathcal{O}(1))^d \cap V\right) = m[\pi(V)]
$$
for some positive integer $m$.
(Note that the dimension of $V$ is $\dim Y-1$
and that $\dim \pi(V) \leq \dim Z_t = \dim Z-1 \leq \dim Y-2$,
so that $d \geq 1$.)
Now suppose $V_1$, $V_2$, \dots, $V_n$
are those components of $D$
for which $d$ achieves its minimum value $d_0$.
Then $s(N_{Z_t}Y_t)$ and $s((N_{Z}Y)_t)$ will differ in dimension
$d_0$ by a \emph{positive} linear combination of
$[\pi(V_1)]$, $[\pi(V_2)]$, \dots, $[\pi(V_n)]$.
Since $Z_t$ is projective, this linear combination
cannot be the zero class.
\end{proof}

\begin{proof}[Proof of Theorem~\ref{mainthm}]
Let $C$ be an elliptic curve.  Then the morphism
$\tilde{Y} := Y \times C \times C \to T$
(projection onto the first factor followed by
$Y \to T$)
is a flat projective morphism from an
equidimensional scheme.
The subscheme $\tilde{Z} := Z \times C$, embedded
in $\tilde{Y}$
via $(z,c) \mapsto (z,c,c)$, is nowhere dense in
$\tilde{Y}$, and its normal cone
$$
N_{\tilde{Z}}\tilde{Y} \cong N_{Z}Y \times N_{C}(C \times C) =
N_{Z}Y \times TC
$$
is internally flat over $T$ since $N_{Z}Y$ is.
The Segre class of $N_{\tilde{Z}}\tilde{Y}$
is obtained from that of $N_{Z}Y$
by pullback via the projection of
$Z \times C \to Z$.  Similarly, $s(N_{\tilde{Z}_t}\tilde{Y}_t)$
is obtained from $s(N_{Z_t}Y_t)$ by pullback via the morphism
$Z_t \times C \to Z_t$.  Therefore, the Segre class
$s(N_{\tilde{Z}}\tilde{Y})$ specializes to $s(N_{\tilde{Z}_t}\tilde{Y}_t)$.

\par
By Lemma~\ref{bootstrap}, we see that $\p N_{\tilde{Z}}\tilde{Y}$ 
is flat over $T$ and
$(\p N_{\tilde{Z}}\tilde{Y})_t = \p N_{\tilde{Z}_t}\tilde{Y}_t$.
Now $\p N_{\tilde{Z}}\tilde{Y}$ is naturally isomorphic to the
product of the projective completion of $N_{Z}Y$ with $C$.
Hence $\p N_{Z}Y \oplus 1$ is flat over $T$ and
$(\p N_{Z}Y)_t = \p N_{Z_t}Y_t$.  Thus $N_{Z}Y$ is flat
over $T$
and $(N_{Z}Y)_t = N_{Z_t}Y_t$.
\end{proof}

\section{Examples:  tangent star cones} \label{examples}

\par
As mentioned in the introduction, the \emph{tangent star cone} 
of a scheme $X$ is the normal cone $N_{X}(X\times X)$ where $X$ is 
embedded in $X\times X$ as the diagonal.  In particular, if $I$ is 
the ideal sheaf of $X$ in $X\times X$, then
$$
TS(X) = \Spec(\oplus_{j \geq 0}I^j/I^{j+1}).
$$
The \emph{projectivized tangent star cone}
$$
\p TS(X) = \Proj(\oplus_{j \geq 0}I^j/I^{j+1}),
$$
is the exceptional divisor of the blowup $Bl_{X}(X\times X)$.  
When $X$ is a reduced subscheme of $\aff^n$, 
we can understand $\p TS(X)$ as the 
closure of the image of the map 
$(X\times X) \setminus X \to \aff^n \times \aff^n \times \p^{n-1}$ defined by
$$
(p_1, p_2) \mapsto (p_1, p_2, \overline{p_1p_2})
$$
where $\overline{p_1p_2}$ denotes the line through the origin 
parallel to the line through the (distinct) points $p_1$ and $p_2$.  
The fiber of $\p TS(X)$ over a point of $X$ consists, 
set-theoretically, of limiting secants, called the \emph{tangent 
star} by K.~Johnson in \cite{johnson}.
For an extensive study of tangent star cones, see \cite{tanstars}.

\par
The tangent star construction above carries over entirely analogously 
to the relative case:  if $X$ is a scheme over a nonsingular variety 
$T$, the \emph{relative tangent star cone}, denoted $TS(X/T)$, is 
$N_{X}(X\times_T X)$.  As above, we have 
$TS(X/T) = \Spec(\oplus_{j \geq 0}I^j/I^{j+1})$ where $I$ denotes 
the ideal sheaf of the diagonal copy of $X$ in $X\times_T X$.

\par
If we apply Theorem~\ref{mainthm} using $X \times_{T} X \to T$
and letting $Z$ be the diagonal copy of $X$, we immediately obtain 
the following result for tangent star cones.

\begin{corollary} \label{tangstarcor}
Suppose that:
\begin{enumerate}
   \item 
   $X \to T$ is a flat, projective morphism from an
   equidimensional scheme to a nonsingular curve germ (with
   closed point $t$);
	 
   \item 
   the tangent star cone $TS(X/T)$ is internally flat over $T$;
   
   \item 
   the Segre class $s(TS(X/T))$ specializes to $s(TS(X_t))$.
\end{enumerate}
\noindent
Then $TS(X/T)$ is flat over $T$, and its fiber over $t$
is  $TS(X_t)$.
\end{corollary}

\par

The Segre classes of relative tangent star cones were studied in
\cite{ky}. We quote Theorem 3 from that paper, using the notation
$s_k$ for the codimension $k$ component of the Segre class.

\begin{theorem} \label{kythm}
Suppose that $Y \to T$ is a smooth morphism to a nonsingular variety.
Suppose that $X$ is a purely codimension $r$ subscheme of $Y$,
and that the composite $X \to T$ is proper and flat. For a closed point
$t$ of $T$, let $X_t$ be the fiber.
If the Segre classes $s_0(TS(X/T))$, \dots, $s_{r-1}(TS(X/T))$ specialize
to the corresponding classes $s_0(TS(X_t))$,\dots, $s_{r-1}(TS(X_t))$,
then the same is true of the higher codimension Segre classes
$s_r(X/T)$, $s_{r+1}(X/T)$, etc.
\end{theorem}

\par
Combining Theorem \ref{kythm} with Corollary~\ref{tangstarcor}, we 
can strengthen that corollary in many instances.

\par
\begin{corollary} \label{improved}
Suppose that:
\begin{enumerate}
   \item 
   $Y \to T$ is a smooth morphism to a nonsingular curve germ;
   
   \item 
   $X$ is a purely codimension $r$ subscheme of $Y$,
   and the composite $X \to T$ is projective and flat;
   
   \item 
   the tangent star cone $TS(X/T)$ is internally flat over $T$;
   
   \item 
   the Segre classes $s_0(TS(X/T)), \dots, s_{r-1}(TS(X/T))$ 
   specialize to the corresponding classes
   $s_0(TS(X_t)),\dots,s_{r-1}(TS(X_t))$.
\end{enumerate}
\noindent
Then $TS(X/T)$ is flat over $T$, and its fiber over $t$
is  $TS(X_t)$. 
\end{corollary}
\par
If $X$ is a hypersurface in $Y$, then there are two especially 
pleasant features. First, by Corollary~\ref{improved}, one
needs to investigate only the top Segre class  $s_0(TS(X/T))$.
If $X$ has irreducible components $X_1,X_2,\dots$ with corresponding
geometric multiplicities $m_1,m_2,\dots$,
then 
$$
s_0(TS(X/T)) = \sum (m_k)^2 [X_k].
$$
(See Theorem~4 of \cite{ky}.)
Thus $s_0(TS(X/T))$ specializes to $s_0(TS(X_t))$ if and only if the
following two criteria are satisfied:
\begin{enumerate}
    \item 
    each component $X_k$ specializes to a reduced subscheme 
    of the fiber $X_t$;
    
    \item 
    the specializations of distinct components have no 
    components in common.
\end{enumerate}
Following the terminology of \cite{gkflat}, we say that
the components of $X$ \emph{do not coalesce} when criteria 1 and 2 
hold.
\par
The second pleasant feature is that $TS(X/T)$ is auto\-matically
Cohen-Macaulay.  Thus we have the following corollary.

\begin{corollary} \label{hypimproved}
Suppose that:
\begin{enumerate}
   \item 
   $Y \to T$ is a smooth morphism to a nonsingular curve germ;
   
   \item 
   $X$ is a hypersurface in $Y$,
   and the composite $X \to T$ is projective and flat;
   
   \item 
   the components of $X$ do not coalesce under specialization. 
\end{enumerate}
\noindent
Then $TS(X/T)$ is flat over $T$, and its fiber over $t$
is  $TS(X_t)$. 
\end{corollary}

\par
In \cite{gkflat} the proof that $TS(X/T)$ is Cohen-Macaulay
is entangled with the proof of Corollary~\ref{hypimproved},
and both proofs employ the
following explicit local description of the tangent star cone.
Suppose $Y$ is a subvariety of $\aff^n \times T$. 
Let $x_1,\ldots,x_n$ be coordinates on $\aff^n$, 
and let $u_1,\ldots,u_m$ be coordinates for 
the tangent bundle with respect to 
$\partial/\partial x_1,\ldots,\partial/\partial x_n$.
Suppose that $X$ is 
defined in $Y$ by the equation $f(x_1,\ldots,x_n,t) = 0$ and write this 
polynomial as
$$
f = \prod_{k=1}^s {f_k}^{r_k}
$$
where the $f_k$'s are reduced, irreducible, and distinct.
Define the \emph{polarization operator} $P$ by
$$
P = \sum_{i=1}^n u_i \frac{\partial}{\partial x_i},
$$
and let $P^d$ denote its $d$th iterate.
Let
$$
S_m f = \left(\prod_{r_k < m}{f_k}^{r_k}\right)^2 P^{2m-1} 
        \left(\prod_{r_k \geq m}{f_k}^{r_k+m-1}\right).  
$$
(Note that $S_1f = Pf$ and $S_m f = 0$ for $m$ sufficiently large.)  
Then the tangent star cone is defined inside the tangent bundle
of $Y$ by the equations
$S_1 f = \ldots = S_m f = 0$
for $m$ sufficiently large.
We would like to find a more conceptual argument that
$TS(X/T)$ is Cohen-Macaulay, thus disentangling the
two proofs and, we hope, giving insight into other situations
where one should expect the relative tangent star cone
to be Cohen-Macaulay (or, more weakly, to satisfy Serre's $S_1$
condition).
\par
In the remainder of this section we look at various other
examples of tangent star cones.
\par

\begin{example} \label{notCMex}
In view of the hypersurface case,
it is tempting to conjecture that if $X \subset Y$ is a 
local complete intersection, then $TS(X)$ is Cohen-Macaulay.  Such a 
conjecture is false, however.  For  example, let $X$ be 
the complete intersection in $\aff^3$ defined by the equations 
$xy = z(z-x) = 0$; then $X$ is a union of three lines, one of 
which is thickened.
Using $a,b,c$ as coordinates with respect to 
$\partial/\partial x, \partial/\partial y, \partial/\partial z$
and calculating using CoCoA \cite{cocoacite}, we find that
the ideal $I$ of the tangent star cone in $R=k[x,y,z,a,b,c]$ is
$$
I = \langle xy, z(z-x), za+xc-2zc, ya+xb, c^2(2xb-zb+yc), 
            bc^2(a^2-2ac+c^2)\rangle.
$$
Then a resolution calculation for $R/I$ shows that its projective dimension
is $5$. But the codimension of $TS(X)$ inside the tangent bundle of $\aff^3$
is $4$. Thus the tangent star cone is not Cohen-Macaulay.
\end{example}

\begin{example} \label{notCMrelex}
Let $X$ be the flat family in $\aff^3$ defined over the affine line 
$T$ by $xy = z(z-tx) = 0$; all fibers except the central one
are isomorphic to the scheme in 
Example~\ref{notCMex}. 
The ideal of the relative tangent star cone
in $R=k[x,y,z,a,b,c,t]$ is
\begin{equation*}
\begin{split}
I = \langle xy, \,& \, z(z-tx), zat+xct-2zc, ya+xb, \\
                & c^2(2xbt-zb+yc), bc^2(a^2t^2-2act+c^2)\rangle.
\end{split}   
\end{equation*}
This family fails to be internally flat over $T$ (and hence 
cannot be Cohen-Macaulay either) by the criterion 
given in Proposition~\ref{intflatprop}:  use 
$$
J = \langle xy, z(z-tx), ya+xb, 
       y^2c^3+a^2bc^2t^2-2abc^3t+bc^4+x^2ct +z^2a-2xzc \rangle
$$
for ``test ideal'' and compare $(I:t)\cap (J:(J:I))$ and $I$.
\end{example}

\begin{example} \label{exampleF}
This example shows that the internal flatness hypothesis cannot be 
omitted from Theorem~\ref{mainthm}.
Let $X$ be the flat family in $\aff^3$ defined over the affine line 
$T$ by the ideal 
$$
\langle x, z\rangle \cap \langle y, z\rangle
            \cap \langle x-y, z-tx\rangle
= \langle z(x-y), xy(x-y), z(z-ty), y(z-tx) \rangle.
$$
The general member of $X$ is the union of three concurrent lines; the 
special member has planar reduction.
The ideal of $TS(X/T)$ is
\begin{equation*}
\begin{split}
I = \langle & z(x-y), xy(x-y), z(z-ty), y(z-tx), \\
         & za - zb + xc - yc, zbt + yct - 2zc, yat+xbt-zb-yc,\\ 
	 & 2xya - y^2a + x^2b - 2xyb, c(abt^2 - act - bct + c^2) \rangle.
\end{split}   
\end{equation*}
This family fails to be internally flat over $T$ by the criterion 
given in Proposition~\ref{intflatprop}:  we use 
\begin{equation*}
\begin{split}
J = \langle & z(x-y), y(z-tx), za-zb+xc-yc, \\
            & (zbt+yct-2zc)+(yat+xbt-zb-yc) +\\ 
	    & (xya-\tfrac{1}{2}y^2a+\tfrac{1}{2}x^2b-xyb)
	       +(abct^2-ac^2t-bc^2t+c^3) \rangle
\end{split}   
\end{equation*}
for test ideal.  
\par
We now show that $s(TS(X/T))$ specializes to $s(TS(X_t))$. 
Writing the general member $X_t$ of the family as the 
union $X_1 \cup X_2 \cup X_3$ of three lines and using 
Theorem 4 of \cite{ky}, we obtain
\begin{equation*}
\begin{split}
s(TS(X_t)) &= s(TS(X_1))+s(TS(X_2))+s(TS(X_3)) \\
             &\quad + 2s(X_1\cap X_2, X_1\times X_2)
	        + 2s(X_1\cap X_3, X_1\times X_3) \\
	     &\quad + 2s(X_1\cap X_3, X_1\times X_3).
\end{split}   
\end{equation*}
Now $X_i\times X_j \cong \aff^2$ and $X_i\cap X_j$ is the origin.  
Hence $s(X_i\cap X_j, X_i\times X_j) = [p]$, the class of a point.  
Therefore, 
\begin{align*}
s(TS(X_t)) &= [X_1]-2[p] + [X_2]-2[p] + [X_3]-2[p] + 6[p] \\
	   &= [X_1] +[X_2] +[X_3] = [X_t].
\end{align*}
Exactly the same calculation applies when $t = 0$, so that 
$s(TS(X_0)) = [X_0]$.  Thus the Segre class specializes.
\end{example}

\end{document}